\documentclass[10pt]{article}

\usepackage{a4wide}
\usepackage{amssymb}
\usepackage{amsfonts}
\usepackage{amsmath}
\input xy
\xyoption{arrow} \xyoption{matrix}

\date{}

\newtheorem{proposition}{Proposition}[section]
\newtheorem{theorem}[proposition]{Theorem}
\newtheorem{lemma}[proposition]{Lemma}
\newtheorem{example}[proposition]{Example}
\newtheorem{definition}[proposition]{Definition}
\newtheorem{corollary}[proposition]{Corollary}

\def\der{\partial }

\def\nFM0{{\nu }_{F,M_0}}
\def\nFN0{{\nu }_{F,N_0}}
\def\nGN0{{\nu }_{G,N_0}}

\def\N0{ {\bf N}_0 }

\def\Xpm{X^{\pm }}

\def\s{\sigma}

\def\l1{{\lambda}_1}

\def\a{\alpha}
\def\a0{ {\alpha }_0}
\def\a1{ {\alpha }_1}

\def\l{\lambda}
\def\o{\omega}

\def\nFGM0{{\nu }_{F,G,M_0}}

\def\nFN0{{\nu}_{F,N_0}}

\def\sm{{\sigma}^m}

\def\sm1{{\sigma}^{-1}}

\def\smtp1{{\sigma}^{-t+1}}

\def\o{\omega }
\def\S1{S^{-1}}

\def\Xpm1{X^{\pm 1}_1}

\def\sPM1{{\sigma }^{\pm 1}}
\def\sMP1{{\sigma }^{\mp 1 }}

\def\di{{\rm d.ind}}

\def\L{\Lambda}

\def\Ytm1{Y^{t-1}}
\def\Yim1{Y^{i-1}}

\def\CK{{\cal K}}

\def\Aut{{\rm Aut}}

\def\dim{{\rm dim }}

\def\ker{ {\rm ker } }

\def\gr{ {\rm gr} }

\def\SL2Z{ {\rm SL}_2({\bf Z}) }

\def\CR{ {\cal R}}

\def\Gp1{ G^{1 , 1 } }
\def\P11{ P^{-1 , 1 } }
\def\Pp1{ P^{1 , 1 } }

\def\nCLsr{{}^\nu\kern-2pt {\cal L}^{\sigma , \rho  }}
\def\nP{{}^\nu \kern-2pt P}
\def\nL{{}^\nu\kern-2pt L}
\def\nLL{{}^\nu\kern-2pt \Lambda}
\def\nPsr{{}^\nu\kern-2pt P^{\sigma , \rho  }}
\def\nLsr{{}^\nu\kern-2pt L^{\sigma , \rho  }}
\def\nuCL{{}^\nu\kern-2pt  {\cal L}}
\def\nCLsr{{}^\nu\kern-2pt {\cal L}^{\sigma , \rho  }}
\def\nCL1m{{}^\nu\kern-2pt {\cal L}^{-1 , 1  }}
\def\x1nu{x^\frac{1}{\nu}}
\def\xm1nu{x^{-\frac{1}{\nu}}}

\def\rad{{\rm rad}}

\def\CR{ {\cal R}}

\def\CB{{\cal B}}

\def\CI{{\cal I}}

\def\CC{ {\cal C}}

\def\nAM0{{\nu }_{{\cal A},M_0}}
\def\nAN0{{\nu }_{{\cal A},N_0}}

\def\CR{ {\cal R }}

\def\bR{\overline{R}}

\def\gn{\mathfrak{n}}
\def\gm{\mathfrak{m}}
\def\gp{\mathfrak{p}}

\def\gr{\mathfrak{r}}

\def\gR{\mathfrak{R}}

\def\SL{{\rm SL}}

\def\Spec{{\rm Spec}}

\def\di!{\frac{\der^i}{i!}}
\def\dik!{\frac{\der^k_i}{k!}}

\def\id{{\rm id}}

\def\N{\mathbb{N}}

\def\0{\overline{0}}
\def\1{\overline{1}}

\def\Ln1{\L_{n,\overline{1}}}

\def\a1{a_{\overline{1}}}

\def\S{\Sigma}

\def\vn1{\overrightarrow{n-1}}

\def\im{{\rm im}}

\def\mJ{\mathbb{J}}
\def\mI{\mathbb{I}}

\def\lann{{\rm l.ann}}

\def\K1{{\rm K}_1}

\def\hmI1{\widehat{\mI_1}}
\def\tmI1{\widetilde{\mI_1}}
\def\tmJ1{\widetilde{\mJ_1}}
\def\hB1{\widehat{B_1}}
\def\hCB1{\widehat{\CB_1}}

\def\rad{{\rm rad}}
\def\s{{\sigma}}

\def\br{\overline{r}}

\def\tor{{\rm tor}}
\def\udim{{\rm udim}}

\def\oG{\overline{G}}

\begin{document}

\author{Volodymyr Bavula  and Vyacheslav Futorny}
\title{Rings of invariants of finite groups when the bad primes exist}

\maketitle

\begin{abstract}

 Let $R$ be a ring (not necessarily with $1$) and $G$ be a finite group of  automorphisms of $R$. The set $\CB(R, G)$ of primes $p$ such that $p\, | \, |G|$ and $R$ is not $p$-torsion free, is called the set of bad primes. When the ring is $|G|$-torsion free, i.e.,  $\CB(R, G)=\emptyset$, the properties of the rings $R$ and $R^G$ are closely connected. The aim of the  paper is to show that this is 
 also true when $\CB(R, G)\neq \emptyset$ under natural conditions on bad primes. In particular, it is shown that the Jacobson radical (resp., the prime radical)  of the ring $R^G$ is equal to the intersection of the Jacobson radical (resp., the prime radical) of $R$ with $R^G$; if the ring $R$ is semiprime then so is $R^G$; if the trace of the ring $R$ is nilpotent then the ring itself is nilpotent; if $R$ is a semiprime ring then $R$ is left Goldie iff the ring $R^G$ is so, and in this case, the ring of $G$-invariants of the  left quotient ring of  $R$ is isomorphic to the left quotient ring of $R^G$ and $\udim (R^G)\leq \udim (R)\leq |G| \udim (R^G)$.


$\noindent $



\end{abstract}


\section{Introduction}
Rings of invariants of finite groups is one of the oldest areas
 of Ring Theory. The main question is how closely properties of the ring and its subring  of invariants are related. Many  classical results  show that the answer is (generally)  affirmative if the order $|G|$ of the group $G$ is a {\em unit} of the ring  or, at least, the ring is $|G|$-{\em torsion free}. 
 Without  these conditions there are many pathological situations.
 
 In particular,  when $|G|$ is a unit of  the ring  the properties of the ring and the subring of invariants are very close, see \cite{Mont}
 \cite{Mont-Pass}, \cite{Lorentz-Mont-Small}, \cite{Jo-Small} and \cite{Alev} for further properties of prime ideals, see also
  \cite{Pascaud} and references therein for some examples. This condition makes it possible to  extend  the classical Noether's result on  affine rings of invariants in commutative rings  to arbitrary 
Noetherian rings $R$ which are  affine over a commutative Noetherian ring $S$ where $G$ is a finite group of $S$-automorphisms of $R$ \cite{Mont-Small}.  In  \cite{Berg},  hereditary and semihereditary properties of the ring and the invariant subring were studies along with the property to be Dedekind among the others. 

A significant step in  study of rings of invariants of finite groups was the classical book of S.Montgomery \cite{Mont-GInv} in which many properties of rings of invariants were studied and plenty of examples were considered especially including many `pathological' examples (when properties of the ring of invariants  differ from the ambient ring). 

Correspondence between  the global dimensions of  the rings and fixed subrings were studied in \cite{Lorentz}.  Integrality of the ring over the invariant subring was discussed in \cite{Pass}.
For representation-theoretical properties of the ring and the invariant subring the interested reader is  referred to \cite{Kraft-Small} and \cite{Dumas}. 

The purpose of this paper is to go beyond of the classical framework  and to establish the conditions under which the invariant subring shares the same  properties as the ambient  ring $R$ in the case when the ring $R$ is {\em not $|G|$-torsion free}, i.e., when the {\em bad primes exist}.

In this paper the following notation is fixed:  $R$ is a ring  not necessarily with $1$, $\Aut(R)$ is its group of (ring) automorphisms, $G$ is a finite subgroup of $\Aut(R)$, $n=|G|$. 
The action of an element $\sigma\in \Aut(R)$ on $R$ is written either as $\sigma(r)$ or $r^{\sigma}$. The subring of $R$, 
$R^G=\{r\in R| \ r^g=r \ \mbox{for all}  \ g\in G\}$, is called the {\em ring of $G$-invariants}. The ideal of $R$, 
$$\tor_n(R):=\{r\in R| \ n^i r=0 \ \mbox{for some} \ i\geq 1\},$$
is called the {\em $n$-torsion ideal} of $R$. The ring $R$ is called {\em $n$-torsion free} if $\tor_n(R)=0$, i.e., the map $n\cdot : R\rightarrow R$, $r\mapsto nr$ is an injection.
 The map $$t=t_G: R\rightarrow R^G, \;\;
r\mapsto \sum_{g\in G}r^g$$ is 
called the {\em trace}.  If $N$ is a normal subgroup of $G$ then $R\supseteq R^N\supseteq (R^N)^{G/N}=R^G$ and we have the map 
$$t_{G/N}: R^N\rightarrow R^G, \, r\mapsto  \sum_{g\in G/N}r^g.$$

The set $\CB(R, G)$  of prime numbers $p$ such that $p| \ |G|$ and $\tor_p(R)\neq 0$ is called the {\em set of bad primes} and the elements of 
$\CB(R, G)$ are called {\em bad primes}.  Clearly, $R$ is $|G|$-torsion free if and only if $\CB(R, G)=\emptyset$.  If the ring  $R$ is $|G|$-torsion free properties of the rings $R$ and $R^G$ 
are closely related, see \cite{Mont-GInv} and below. In this paper, we explore the case when $\CB(R, G)\neq \emptyset$ and find natural conditions under which properties of $R$ and $R^G$ are closely related in a similar way as in the $|G|$-torsion free  case.\\

{\bf When $t_G(R)$ is nilpotent implies $R$ is nilpotent.}
There are two classical results in this direction: [Theorem 1.4, \cite{Mont-GInv} and [Theorem 1.7, \cite{Mont-GInv}], see below. 

For a finite group $G$, set
$$h(G)=\prod_{i=1}^{|G|}\Big( {|G|\choose i}+1 \Big).
$$
The next result is due to G. Bergman and I. M. Isaacs \cite{11}.

\begin{theorem}[Theorem 1.4, \cite{Mont-GInv}]\label{Th1.4Mont}
Let $R$ be a ring and $G$ be a finite subgroup of $\Aut(R)$. If $R$ is  $|G|$-torsion free and $t(R)^d=0$ for some $d\geq 1$ then $R^{h(G)^d}=0$. 
\end{theorem}

Let $p$ be a prime number such that $p| \ |G|$. Then $|G|=p^s m$ for some $s\geq 1$ and $m$ which is relatively prime to $p$. A subgroup $H$ of $G$ is called a {\em $p$-complement} if $|H|=m$. If, in addition, the group $H$ is normal then it is called a {\em $p$-normal complement}. If a $p$-normal complement exists then it contains precisely all the elements of the group $G$ of order not divisible by $p$. Therefore, it is {\em unique} and  denoted   by $N(p)$.

\begin{theorem}[Theorem 1.7, \cite{Mont-GInv}]\label{Th1.7Mont}
Let $R$ be a ring and $G$ be a finite subgroup of $\Aut(R)$  such that $R^G=\{0\}$. Suppose that for every prime $p\in \CB(R, G)$, $G$ has a  $p$-normal complement. Then $R$ is nilpotent. 
\end{theorem}

Theorem \ref{Th1.7Mont} is a  special case of  Theorem \ref{N-1} as  conditions  1-3 of Theorem \ref{N-1} are automatically hold when $R^G=0$.

\begin{theorem}\label{N-1}
Let $R$ be a ring, $G$  a finite subgroup of $\Aut(R)$ and $n=|G|$. Suppose that $\CB(R, G)\neq \emptyset$ and for each $p\in\CB(R, G)$:
\begin{enumerate}
\item  The group $G$ has a  $p$-normal complement $N(p)$;
\item   The ring $R^G$ is $p$-torsion free;
\item   $t_{G(p)}(R^{N(p)})^{d(p)}=0$ for some natural number $d(p)\geq 1$ where $G(p):=G/N(p)$ and $$t_{G(p)}: R^{N(p)} \rightarrow (R^{N(p)})^{G(p)}=R^G, \,  r\mapsto \sum_{\sigma\in G(p)}r^{\sigma}.$$
\end{enumerate}
Then
\begin{enumerate}
\item $R^l=0$ where $l=\max_{p\in \CB(R, G)} \{ l(p)\}$ and  $l(p):=h(N(p))^{h(G(p))^{d(p)}}$. In particular, the ring $R$ is  nilpotent. 
\item For every prime $p\in \CB(R, G)$, the ring $R^{N(p)}$ is a $p$-torsion free, nilpotent ring. Furthermore, $(R^{N(p)})^{m(p)}=0$ where $m(p)=h(G(p))^{d(p)}$.
\end{enumerate}
\end{theorem}

The proof of Theorem \ref{N-1} is given in Section \ref{Proofs}. 

\

{\bf When $R$ is semiprime implies $R^G$ is semiprime.} We recall the following classical result.

\begin{theorem}[Corrolary 1.5, \cite{Mont-GInv}]\label{C1.5Mont}
Let $R$ be a semiprime ring and $G$ be   a finite subgroup of $\Aut(R)$. If $R$ is $|G|$-torsion free then 
\begin{enumerate}
\item The ring $R^G$ is a semiprime ring. 
\item $t(I)\neq 0$ for all nonzero $G$-invariant left or right ideals $I$ of $R$. 
\end{enumerate}
\end{theorem}

For any ring $R$ if the group $G$ satisfies both properties $1$ and $2$ of Theorem \ref{C1.5Mont} then we say that the group has {\em non-degenerate trace}.

 Theorem \ref{N-2} is an extension  of Theorem \ref{C1.5Mont} to the case when $\CB(R, G)\neq \emptyset$.

\begin{theorem}\label{N-2}
Let $R$ be a semiprime ring and $G$ be   a finite subgroup of $\Aut(R)$. Suppose that $\CB(R, G)\neq \emptyset$ and for every $p\in\CB(R, G)$:
\begin{enumerate}
\item The group $G$ has a  $p$-normal complement $N(p)$, and 
\item The ring $R^G$ is $|G|$-torsion free.
\end{enumerate}
Then
\begin{enumerate}
\item The ring $R^G$ is a semiprime ring.
\item  For all nonzero $G$-invariant left or right ideals $I$ of $R$, $t(I)\neq 0$. Furthermore, $t(I)^i\neq 0$ for all $i\geq 1$.
\end{enumerate}
\end{theorem}

For a ring $R$, let $\CC_R$ be the set of regular elements of $R$ ($\CC_R$ is the set of non-zero-divisors of $R$). The rings $Q_l(R):=\CC_R^{-1}R$ and $Q_r(R)=R\CC_R^{-1}$ are called the {\em (classical) left and right quotient rings} of $R$, respectively. A ring $R$ has {\em finite left uniform dimension}, $\udim (R)<\infty$, if it does not contain an  infinite direct sum  of nonzero left ideals. A ring $R$ is called a {\em left Goldie ring} if it has finite left uniform dimension and satisfies the {\em ascending chain condition} for left annihilator ideals. Recall that for a  non-empty subset $X$ of $R$,  the left ideal of $R$,   $\lann_R(X):=\{ r\in R\, | \, rX=0\}$, is called the {\em left annihilator} of the set $X$ and a left ideal of this kind is called a {\em left annihilator ideal} of $R$. 

\begin{corollary}\label{a8Apr17}
Let $R$ be a semiprime ring and $G$ be a finite subgroup of $\Aut (R)$. Suppose that $\CB (R,G)\neq 0$ and conditions 1 and 2 of Theorem \ref{N-2} hold. Then 
\begin{enumerate}
\item The ring $R$ is a left (resp., right) Goldie ring iff the ring $R^G$ is so. In this situation, $Q_l(R)^G=Q_l(R^G)$ (resp.,  
$Q_r(R)^G=Q_r(R^G)$), $\CC_{R^G}\subseteq \CC_R$ and $Q_l(R)=\CC_{R^G}^{-1}R$ (resp., $Q_r(R)=R\CC_{R^G}^{-1}$). 
\item Let $\udim$ be either a left or right uniform dimension. Then $\udim (R)<\infty $ iff $\udim (R^G)<\infty$, and in this case, $$\udim (R^G)\leq \udim (R)\leq |G|\udim (R^G).$$
\end{enumerate}

\end{corollary}

{\bf The equality $\gn(R^G)=R^G\cap \gn(R)$.} For a ring $R$, we denote $\gn(R)$ its {\em prime radical}, i.e., $\gn(R):=\cap_{\gp\in \Spec(R)}\gp$ where $\Spec (R)$ is the prime spectrum of the ring $R$.  The next theorem appears in \cite{Mont-GInv} with similar results in \cite{26} and \cite{43}.

\begin{theorem}[Theorem 1.9, \cite{Mont-GInv}]\label{Th1.9Mont}
Let $R$ be $|G|$-torsion free. Then $\gn(R^G)=R^G\cap \gn(R)$. 
\end{theorem}

The next theorem shows that the same result holds when $\CB(R, G)\neq \emptyset$ under the assumptions of Theorem \ref{N-2} but for the ring $R/\gn (R) $.

\begin{theorem}\label{4Apr17}
Let $R$ be a ring, $G$ be  a finite subgroup of $\Aut(R)$, $\bR=R/\gn(R)$ and $\oG$ be the image of the group $G$ under the  group homomorphism $G\rightarrow \Aut(\bR)$, $g\mapsto \bar{g}$ where 
$\bar{g}(r+\gn(R)):=g(r)+\gn(R)$. Suppose that either $\CB(\bR,\oG )= \emptyset$    or $\CB(\bR, \oG)\neq \emptyset$ and for every $p\in \CB(\bR, \oG)$:
\begin{enumerate}
\item The group $\oG$ has a  $p$-normal complement $N(p)$, and 
\item The ring $\bR^{\oG}$ is $|G|$-torsion free.
\end{enumerate}
Then 
\begin{enumerate}
\item $\gn(R^G)=R^G\cap \gn(R)$.
\item The rings  $\bR^{\overline{G}}$  and $\overline{R^G}$ are semiprime and $|G|\bR^{\overline{ G}}\subseteq \overline{R^G}\subseteq \bR^{\oG}.$
\end{enumerate}

\end{theorem}

The proof of Theorem \ref{4Apr17} is given in Section 2.

\

{\bf The radical of the ring of invariants.} 
When $|G|^{-1}\in R$, the (Jacobson) radicals $\rad (R)$ and $\rad(R^G)$ of the rings $R$ and $R^G$ are closely related, see Theorem \ref{Th1.4-Mont} 
below which is due to S.~Montgomery.

\begin{theorem}[Theorem 1.4, \cite{Mont-GInv}]\label{Th1.4-Mont}
Let $R$ be a  ring and $G$ be   a finite subgroup of $\Aut(R)$. Suppose that $|G|^{-1}\in R$. Then 
$$\rad(R^G)=\rad(R)\cap R^G.$$ 
\end{theorem}

 In general, even for domains  the condition  `$|G|^{-1}\in R$' in Theorem \ref{Th1.4-Mont} cannot be replaced by the weaker  condition that `the ring $R$ is $|G|$-torsion free', \cite{53} (see also \cite{78} for a simpler example).  The theorem below shows that  under certain conditions the same result as in Theorem \ref{Th1.4-Mont} holds when the set of bad primes is a nonempty set, see Section \ref{RAD-INV} for details. 

\begin{theorem}\label{B5Apr17}
Let $R$ be a  ring, $\pi: R\rightarrow \bR:= R/\rad(R)$, $r\mapsto \br:=r+\rad(R)$. Let $G$ be   a finite subgroup of $\Aut(R)$ and $\overline G$ be its image under the group homomorphism 
$\Aut(R)\rightarrow \Aut(\bR)$, $g\mapsto \overline g$ where ${\overline g}(\br)=\overline{g(r)}$.  Suppose that $\overline{R^G}=\bR^{\overline G}$, the group $\overline G$ is either a left or right proper splitting group for the ring $\bR$ and either
\begin{enumerate}
\item The ring $\bR$ is $|\overline G|$-torsion free, or
\item $\CB(\bR, \overline G)\neq \emptyset$ and conditions 1 and 2 of Theorem \ref{N-2} hold for the ring $\bR$ and the group $\overline G$. 
\end{enumerate}
Then $\rad(R^G)=\rad(R)\cap R^G$.
\end{theorem}

{\bf The ring $R^G$ is a semisimple Artinian ring.} The next theorem is an  old result of Levitzki, \cite{46}, 

\begin{theorem}\label{Th-Levitzki}
Let $R$ be a  ring and $G$ be   a finite subgroup of $\Aut(R)$. Suppose that $|G|^{-1}\in R$ and the ring $R$ is a semisimple Artinian ring. Then the ring $R^G$ is  a semisimple Artinian ring.
\end{theorem}
The theorem below shows that under certain conditions the same result is true  in case the set of bad primes is a non-empty set. 

\begin{theorem}\label{8Apr17}
Let $R$ be a semisimple Artinian  ring and $G$ be a finite subgroup of $\Aut(R)$.
  Suppose that  the group $G$ is either a left or right proper splitting group for the ring $R$ and either
\begin{enumerate}
\item The ring $R$ is $|G|$-torsion free, or
\item $\CB(R, G)\neq \emptyset$ and conditions 1 and 2 of Theorem \ref{N-2} hold for the ring $R$. 
\end{enumerate}
Then $R^G$ is  a semisimple Artinian  ring.
\end{theorem}

The proof of Theorem \ref{8Apr17} is given in Section \ref{RAD-INV}.

\begin{corollary}\label{b8Apr17}
Let $R$ be a semiprime ring and $G$ be a finite subgroup of $\Aut (R)$. Suppose that $\CB (R,G)\neq 0$ and conditions 1 and 2 of Theorem \ref{N-2} hold. Then the ring 
$R$ is a semisimple Artinian ring iff  the ring $R^G$ is so, and in this case, $$\udim (R^G)\leq \udim (R)\leq |G|\udim (R^G).$$
\end{corollary}

{\bf Proof.} The corollary follows from Corollary \ref{a8Apr17} and the fact that semisimple Artinian ring coincides with its (left and right) quotient ring.  $\Box$


\section{Proofs  of Theorem \ref{N-1}, Theorem \ref{N-2} and Theorem \ref{4Apr17}}\label{Proofs}

The aim of this section is to give proofs of Theorem \ref{N-1}, Theorem \ref{N-2} and Theorem \ref{4Apr17}. 

\

Let $p$ be a prime number. A finite group $H$ is called a $p$-{group} if $|H|=p^s$ for some natural number $s$. \\

{\bf Proof of Theorem \ref{N-1}.}  Let $p\in\CB(R, G)$. Then $n=p^s m$ for some natural number  $m$ not divisible by $p$. Then $|N(p)|=m$ and the factor group $G(p)$ is a finite $p$-group since $|G(p)|=p^s$, $s\geq 1$.

 {\em (i) The ring $R^{N(p)}$ is $p$-torsion free:}
The group $G(p)$ acts on the ring $R^{N(p)}$. Suppose that $T(p):= \tor_p(R^{N(p)})\neq 0$. We seek a contradiction.  Then the group $G(p)$ acts on  $T(p)$ and also on its nonzero abelian subgroup $\CK(p):=\ker_{T(p)}(p\cdot)$, which is an $\mathbb F_p$-module where $\mathbb F_p=\mathbb Z/p\mathbb Z$ is the finite field that contains $p$ elements.  Fix a nonzero element $r$ of $\CK(p)$. Then the 
$\mathbb F_p$-module $V=\mathbb F_p G(p) r$ is an abelian $p$-group (since $|V|=p^{\dim_{\mathbb F_p}(V)}$). Then (it is well-known that) the action of the finite $p$-group $G(p)$ on a finite abelian $p$-group $V$ 
has a {\em nonzero} fixed point, say $v(p)$. Since $v(p)\in (R^{N(p)})^{G(p)}=R^G$, we have $v(p)\in R^G\cap \CK(p)=0$ (by assumption 2),  a contradiction.

 {\em (ii)   $\Big(|G(p)|R^{N(p)}\Big)^{m(p)}=0$ where $m(p)=h(G(p))^{d(p)}$:}  Applying [Proposition 1.3, \cite{Mont-GInv}] to the group $G(p)$ and the ring $R^{N(p)}$, we have the inclusion 
$$\Big(|G(p)|R^{N(p)}\Big)^{h(G(p))^{d(p)}}\subseteq |G(p)|(R^{N(p)})' t_{G(p)}(R^{N(p)})^{d(p)}R^{N(p)}=0$$ (by assumption 3). Here $(R^{N(p)})'$ denotes a ring that is obtained from the ring $R^{N(p)}$ by adding 1. 

\

 {\em (iii)   $\Big(R^{N(p)}\Big)^{m(p)}=0$:} Follows from the statement {\em (ii)} and assumption 2. \\
 
 \
 Let $l(p)=h(N(p))^{m(p)}$. Then 

\
 
  {\em (iv) $(|N(p)|R)^{l(p)}=0$:} Applying [Proposition 1.3, \cite{Mont-GInv}]  to the group $N(p)$ and the ring $R$ we have that 
  $$\Big(|N(p)|R\Big)^{h(N(p))^{m(p)}}\subseteq |N(p)|R' t_{N(p)}(R)^{m(p)}R\subseteq |N(p)|R' \Big(R^{N(p)}\Big)^{m(p)}R=0,$$ 
  by the statement {\em (iii)} where  $R'$   denotes a ring that is obtained from the ring $R$ by adding 1.

\

{\em (v)   $R^{l(p)}$ is $p$-torsion free}, by the statement {\em (iv)} and the fact that $p$ does not divide $|N(p)|$. \\
\

{\em (vi)   $R^{l}$ is $n$-torsion free}, by the statement {\em (v)} and since $R^l\subseteq R^{l(p)}$ for all $p\in\CB(R, G)$.

\

{\em (vii)   $R^{l}=0$}, by the statements {\em (iv)} and  {\em (vi)} and since $R^l\subseteq R^{l(p)}$. \;\; $\Box$ 

\

{\bf Proof of Theorem \ref{N-2}.}  1. Suppose that 
 $J$ be a nonzero nilpotent ideal $R^G$. We seek a contradiction. Fix a nonzero element $a$ of 
$J$. Then $(aR^G)^m=0$ for some $m\geq 1$, the right ideal $aR$ of $R$ is $G$-invariant (for all $g\in G$, $g(aR)=g(a)g(R)=aR$) and $t(aR)=at(R)$ is nilpotent, since $(at(R))^m\subseteq (aR^G)^m=0$, i.e., $t(aR)^m=0$. 

The ring $\gR=aR$ satisfies conditions {\em 1-3} of Theorem \ref{N-1}. Indeed, conditions {\em 1-2} of Theorem \ref{N-1} follow from  conditions {\em 1-2} of the theorem. 
 Let $p\in\CB(R, G)$. Then 
 $$0=t(aR)^m\supseteq t((aR)^{N(p)})^m=t(\gR^{N(p)})^m=(|N(p)|t_{G(p)}(\gR^{N(p)}))^m.$$
 So, $T(p):=t_{G(p)}(\gR^{N(p)})^m$ is $|N(p)|$-torsion, hence $|G|$-torsion.  Since
 $$t_{G(p)}(\gR^{N(p)}))^m\subseteq t_{G(p)}(R^{N(p)}))^m\subseteq (R^G)^m,$$  $T(p)=0$ by assumption {\em 2}, and so condition {\em 3} of Theorem \ref{N-1} holds.
Therefore, by Theorem \ref{N-1}, the ring (the right ideal of $R$) $\gR=aR$ is nilpotent which is a contradiction ($R$ is semiprime). 

2. Suppose that   $t(I)^i=0$ for some nonzero $G$-invariant left or right ideal $I$ of $R$ and  $i\geq 1$. We seek a contradiction. The ring $I$ satisfies conditions 
{\em 1-3} of Theorem \ref{N-1}: conditions 
{\em 1,2} of Theorem \ref{N-1} follow from the conditions 
{\em 1,2} of the theorem.  Let $p\in\CB(R, G)$. Then 
$$0=t(I)^i\supseteq t(I^{N(p)})^i=(|N(p)|t_{G(p)}(I^{N(p)}))^i.$$
 So, $S(p):=t_{G(p)}(I^{N(p)})^i$ is $|N(p)|$-torsion, hence $|G|$-torsion.  Since $S(p)\subseteq R^G$,  we must have $S(p)=0$, by assumption {\em 2}, and so 
 condition {\em 3} of Theorem \ref{N-1} holds.  Now by Theorem \ref{N-1}, the ring (the left or right ideal of $R$) $I$ is nilpotent, which is a contradiction ($R$ is semiprime). $\Box$

\

{\bf Proof of Corollary \ref{a8Apr17}.} Corollary \ref{a8Apr17} follows from Theorem \ref{N-2} and [Theorem 5.3, \cite{Mont-GInv}].  $\Box$

\

For a ring $R$, we denote by $\CI(R)$ the set of ideals of $R$. The inclusion of rings $R^G\subseteq R$ yields the {\em restriction} and {\em extension} maps,
$$\CI(R)\rightarrow \CI(R^G), \ I\mapsto I^r:= R\cap I\;\;{\rm
and}\;\;\CI(R^G)\rightarrow \CI(R), \ J\mapsto J^e:= (J)= R J R.$$\\

{\bf Proof of Theorem \ref{4Apr17}.} The set $\CR:=\{I\in \CI(R) | \  I\cap R^G\subseteq \gn(R^G)\}$ is a non-empty set as $\{0\}\in \CR$. By Zorn's Lemma, the set $\max(\CR)$ of maximal (with respect to $\subseteq$) elements of $\CR$ is a non-empty set.

{\em (i) All elements of $\max(\CR)$ are semiprime ideals of $R$:} Let $I\in \max (\CR )$. We have to show that if  $J^2\subseteq I$  for some ideal $J$ of  $R$ containing $I$ then $J=I$. Since
$(J\cap R^G)^2\subseteq J^2\cap R^G\subseteq I\cap R^G\subseteq \gn (R^G)$, we must have $J\cap R^G\subseteq \gn (R^G)$, and so $J\in \CR$. Now, $J=I$, by the maximality of $I$.

{\em (ii) $R^G\cap \gn(G)\subseteq \gn(R^G)$:} Take $\gm\in \max (\CR)$. Since the ideal $\gm$ of $R$ is semiprime (the statement {\em (i)}),  $\gn(R)\subseteq \gm $. Hence,  
$R^G\cap \gn(G)\subseteq R^G\cap \gm\subseteq \gn(R^G)$ (since $\gm\in \CR$).

{\em (iii) The ring $\bR^{\oG}$ is semiprime:} If $\CB(\bR, \oG)= \emptyset$, i.e., the ring $\bR$ is a semiprime, $|G|$-torsion free ring, then $\bR^{\oG}$ is a semiprime ring, by Theorem \ref{C1.5Mont}. 
 If $\CB(\bR, \oG)\neq \emptyset$ then the ring $\bR^{\oG}$ is a semiprime ring, by Theorem \ref{N-2}.

 {\em (iv)  $|G|\bR^{\oG}\subseteq \overline{R^G}\subseteq \bR^{\oG}$:} The inclusions are obvious.
 
 {\em (v) The ring $\overline{R^G}\simeq R^G/R^G\cap \gn(G)$ is semiprime:} Let $J$ be a nonzero nilpotent ideal of the ring  $\overline{R^G}$. Then the set $$\tilde{J}:=\{a\in \bR^{\oG}\, | \ |G|^i a\in J \ \mbox{for some} \; i\geq 1\}$$ 
 is an ideal of the ring  $\bR^{\oG}$ such that 
 $$0\neq |G|J\subseteq  |G|\tilde{J}\subseteq |G|\bR^{\oG}\subseteq \overline{R^G} \subseteq \overline {R^G},$$ by the statement {\em (iv)} and assumption {\em 2}.  The  nonzero ideal $ |G|\tilde{J}$ of 
 $\bR^{\oG}$ is a nilpotent ideal (since $\overline{ R^G}\subseteq \bR^{\oG}$ and $J$ is a nilpotent ideal of $\overline {R^G}$). This is contradiction to the statement {\em (iii)}.
 
 {\em (vi) $\gn(R^G)\subseteq R^G\cap \gn(G)$:} The inclusion follows from the statement {\em (v)}. $\Box$

\

\


\section {The radical of the ring of invariants} \label{RAD-INV}
\

The aim of this section is to give proofs of Theorem \ref{B5Apr17} and Theorem \ref{A5Apr17} that connect the (Jacobson) radical of a ring and the radical of the ring of invariants. 

\

{\bf Splitting subrings  and splitting groups.}

\begin{definition}
Let $A\subseteq B$ be rings. The subring $A$ of $B$ is called a {\em splitting subring for $B$} if
$B=A\oplus A'$ for some $A$-subbimodule $A'$ of $B$. The $A$-subbimodule $A'$ is called a {\em splitting $A$-subbimodule}.  
In the particular case when $A=R^G\subset B=R$ where $R$ is a ring and $G$ is a subgroup of $\Aut(R)$, we say that $R^G$ is a {\em splitting subring of $R$} and $G$ is a {\em splitting group} (of automorphisms).  
\end{definition}

In general, splitting $A$-subbimodule is not unique. Indeed, let $A$ be a field and $B$ be an $A$-algebra such that $A\neq B$. Then every $A$-subspace $A''$ of $B$ such that $B=A\oplus A''$ 
is a splitting $A$-subbimodule.

In general situation, suppose that $A$ is a  splitting subring of $B$ and $A'$ is a splitting $A$-subbimodule of $B$.  Let $\Aut_A(B)$  and $\Aut^A(B)$  be the subgroups of $\Aut (B)$ 
 that contains automorphisms $\s\in \Aut(B)$ such that $\s(A)=A$ and $\s_A=\id_A$ in the first case. So,  $\Aut_A(B)$  is a normal subgroup of  $\Aut^A(B)$. For every $\tau\in \Aut^A(B)$, 
 $B=\tau(B)= \tau(A)\oplus \tau(A')=A\oplus \tau(A')$ is a direct of $A$-bimodules. So, $\tau(A')$ a splitting $A$-subbimodule.
 
 The sets $C_B(A)=\{b\in B| \ ba=ab \ \mbox{for all } \ a\in A\}$ and  $N_B(A)=\{b\in B| \ bA=Ab \}$ are called the {\em centralizer} and the {\em normalizer} of $A$ in $B$, respectively. Clearly, 
$C_B(A)\subseteq N_B(A)$ are subrings of $B$. Let $B^{\times}$ be the group of units of $B$. Each unit $u\in B$ determines the inner automorphism $\o_u$ of the ring $B$ given by the rule 
$\o_u(b)=ubu^{-1}$. 

Let $C_B(A)^{\ast}=C_B(A)\cap B^{\times}$ and $N_B(A)^{\ast}=N_B(A)\cap B^{\times}$. Then  $C_B(A)^{\ast}$ is a normal subgroup of $N_B(A)^{\ast}$. Let $c\in C_B(A)^{\ast}$ and 
$\nu\in N_B(A)^{\ast}$. Then $\o_c\in \Aut_A(B)$ and $\o_{\nu}\in \Aut^A(B)$, and so  $B=A\oplus \o_c(A')=A\oplus \o_{\nu}(A')$ where $ \o_c(A')$ and $ \o_{\nu}(A')$ are splitting 
$A$-subbimodules of $B$. 

\begin{example}
Let $R$ be a ring and $G$ be a finite subgroup of $\Aut (G)$ such  that $|G|^{-1}\in R$. Then $R=R^G\oplus B$ is a direct sum of $R^G$-bimodules where $B=(1-e)(R)$ and 
$$e=\frac{1}{|G|}\sum_{g\in G}g$$
(since $e^2=e$, $e(R)=R^G$ and the map $e\cdot : R\rightarrow R$, $r\mapsto e(r)=\frac{1}{|G|}\sum_{g\in G}g(r)$ is a homomorphism of $R^G$-bimodules).
\end{example}

\begin{example}
In the previous example, the condition  $|G|^{-1}\in R$ is a sufficient but not necessary condition for the group $G$ to be splitting. Indeed, take any ring $R$ with $R^G=\{0\}$ and $|G|R=0$, see 
[Example 1.1, p.6, \cite{Mont-GInv}], and any ring $S$ with $|G|S=0$. Let $A=S\times R$ be a direct product of rings. Extend the action of the group $G$ from $R$ to $A$ by the rule 
$s^g=s$ for all $s\in S$ and $g\in G$. Then $A^G=S^G\times R^G=S^G=S $ and $A=A^G\times R$ is a direct sum of $A^G$-bimodules but $|G|A=0$. 

\end{example}

For a ring $R$, let $\CI_l(R)$ and $\CI_r(R)$ be the sets of left and right ideals of $R$, respectively. The maps $\CI_l(R^G)\rightarrow \CI_l(R)$, $J\mapsto J^e : = RJ$ and 
$\CI_r(R^G)\rightarrow \CI_r(R)$, $J\mapsto J^e: = JR$ are called the {\em extension } maps, and the maps $\CI_l(R)\rightarrow \CI_l(R^G)$, $I\mapsto I^r : = I\cap R^G$ and 
$\CI_r(R)\rightarrow \CI_r(R^G)$, $I\mapsto I^r: = I\cap R^G$ are called the {\em restriction} maps. Clearly, $J^{er}\supseteq J$ and $I^{re}\subseteq I$. 

For a left (right) $R$-module $M$, $l_R(M)$ denotes its {\em length}.

\begin{lemma}\label{b6Apr17}
Let $R$ be a ring and $G$ be a finite subgroup of $\Aut (G)$. Suppose that $G$ is a splitting group and $R= R^G \oplus B$ is a direct sum of $R^G$-bimodules. Then 
\begin{enumerate}
\item For all left (right) ideals $J$ of $R^G$, $J^{er}=J$. So, the extension map $J\mapsto J^e$ is an injection. 
\item For all left (right) ideals $J$ of $R^G$, $l_{R^G}(R^G/J)\leq l_R(R/J^e)$.
\end{enumerate}
\end{lemma}

\begin{definition}
A direct sum of $R^G$-bimodules $R=R^G\oplus B$ is called a {\em left (respectively, right) proper splitting} if $e(I)\subseteq I\cap R^G$ for all left (respectively,  right) $G$-invariant ideals $I$ of $R$ where $e$ is the projection of $R$ onto $R^G$.  A group $G$ is called a {\em left (respectively, right) proper splitting group} if there is at least one  left (respectively,  right) proper splitting.
\end{definition}

\begin{example}
Suppose that  $|G|^{-1}\in R$ and $e=\frac{1}{|G|}\sum_{g\in G}g$. Then $R=R^G\oplus B$, where $B=\im(1-e)$,  is a left and right proper splitting and the group $G$ is a left and right proper splitting group.
\end{example}

\begin{lemma}\label{a6Apr17}
Let  $R=R^G\oplus B$  be a direct sum of $R^G$-bimodules and $e$ be the projection onto  $R^G$.  Then $R=R^G\oplus B$ is a left (respectively, right) proper splitting iff for all left (respectively, right) 
$G$-invariant ideals $I$ of $R$, $e(I)=I\cap R^G$. 
\end{lemma}

{\bf Proof. } $(\Rightarrow)$ For all left (respectively, right) 
$G$-invariant ideals $I$ of $R$, $I\cap R^G\subseteq e(I)\subseteq I\cap R^G$, and so $ e(I)= I\cap R^G$.

$(\Leftarrow)$  This implication is obvious.  $\Box$

\

Let $R$ be a ring and $I$ be a left (respectively, right) ideal of $R$. Let $\CI_l(R, I)$ (respectively, $\CI_r(R, I)$ be the set of all left (respectively, right) ideals of $R$ that contain $I$. 

\begin{lemma}\label{c6Apr17}
Let $R$ be a ring and $G$ be a finite subgroup of $\Aut (G)$. Suppose that $R=R^G\oplus B$ is a left (respectively, right) proper splitting and $e$ is the projection of $R$ onto  $R^G$. Then 
\begin{enumerate}
\item For every left (respectively, right) $G$-invariant ideal $I$ of $R$, $I=I\cap R^G$. 
\item For every left (respectively, right) $G$-invariant ideal $I$ of $R$,  the map $\CI_l(R^G, I\cap R^G)\rightarrow \CI_l(R, I)$, $J\mapsto  RJ+I$ (respectively,  $\CI_r(R^G, I\cap R^G)\rightarrow 
\CI_r(R, I)$, $J\mapsto JR+I$) is injective since $(RJ+I)\cap R^G=J$ (respectively,  $(JR+I)\cap R^G=J$).
\item For every left (respectively, right) $G$-invariant ideal $I$ of $R$,  $l_{R^G}(R^G/R^G\cap I)\leq l_R(R^G/I)$.
\item If $R/I$ is an Artinian or Noetherian left (resp., right) $R$-module then the left (resp., right) $R^G$-module $R^G/R^G\cap I$ is so. 
\end{enumerate}
\end{lemma}

{\bf Proof. } 1. By Lemma \ref{a6Apr17}, $e(I)=I\cap R^G\subseteq I$. Notice that for all $i\in I$, $i=e(i)+(1-e)(i)$.  Therefore, $(1-e)(i)=i-e(i)\in I\cap B$ for all $i\in I$, and statement 1 follows.

2.  For every left (respectively, right) $G$-invariant ideal $I$ of $R$, the left (respectively, right) ideal $RJ+I$ (respectively, $JR+I$) of $R$ is $G$-invariant,   
$$e(RJ+I)=e(RJ)+e(I)= e(J\oplus BJ)+ I\cap R^G=J+I\cap R^G=J$$ 
(respectively, $e(JR+I)=J+I\cap R^G=J$).  Now, statement 2 is obvious.

3 and 4.  Statements 3 and 4 follow at once from statement 2.  $\Box$

When $|G|^{-1}\in R$ we have the well-known result, see [Section 1, \cite{Mont-GInv}], in particular [Corollary 1.12, \cite{Mont-GInv}]. 

\begin{corollary}\label{c8Apr17}
Let $R$ be a ring and $G$ be a finite subgroup of $\Aut (G)$. Suppose that $|G|^{-1}\in R$. Then Lemma \ref{c6Apr17}  holds 
 where $e=|G|^{-1}\sum_{g\in G}g$. 
\end{corollary}

So, when $|G|^{-1}\in R$ the ring $R$ is a left or right Artinian or Noetherian iff the ring $R^G$ is so. In general, this result is false if the ring $R$ is not $|G|$-torsion free, see \cite{63,14}.

\begin{theorem}\label{A5Apr17}
Let $R$ be a ring and $G$ be a finite subgroup of $\Aut (G)$. Suppose that $\rad (R)=0$ and the group $G$ is either a left or right proper splitting group and either 
\begin{enumerate} 
\item The ring $R$ is $|G|$-torsion free, or
\item 
$\CB(R, G)\neq \emptyset$ and conditions 1 and 2 of Theorem \ref{N-2} hold. 
\end{enumerate}
Then $\rad(R^G)=0$. 
\end{theorem}

{\bf Proof.}
{\em (i) $R$ is semiprime ring } (since $\gn (R)\subseteq \rad(R)=0$).

{\em (ii) $R^G$ is a semiprime ring:} Suppose that condition 1 (respectively, 2) holds. Then, by Theorem \ref{Th1.4Mont} (respectively, Theorem \ref{N-2}) the ring $R^G$ is a semiprime ring. 

Let us assume that the group $G$ is a left proper splitting group (the `right' case can be dealt in a similar fashion).

{\em (iii) $\rad(R^G)=0$:} In view of the statement {\em (ii)},  it is enough to show that $\gr:=\rad (R^G)$ is nilpotent.  Let $I$ be a maximal regular left ideal of the ring $R$. The ideal $I^{\circ}=\bigcap_{g\in G}g(I)$ is a left $G$-stable  ideal of $R$ ($GI=I$). There is an obvious $R$-module homomorphism 
$$R/I^{\circ}\rightarrow \prod_{g\in G}R/g(I), \, r+ I^{\circ}\mapsto (r+g(I))_{g\in G}.$$
Notice that $V:= R^G/I^{\circ}\cap R^G$ is a  left  $R^G$-submodule of $R/I^{\circ}$. Recall that $G$ is a left proper  splitting group. So, by Lemma \ref{c6Apr17},
$$l_{R^G}(V)\leq l_R(R/I^{\circ})\leq l_R\Big( \prod_{g\in G}R/g(I)\Big)=n:=|G|< \infty.$$

Hence, $\gr^n V=0$, i.e.,  $\gr^n R^G \subseteq I^{\circ}\cap R^G \subseteq I$, and so $\gr^n R^G \subseteq \bigcap I = \rad(R)=0$. Therefore, $\gr$ is a nilpotent ideal of $R^G$ as required. 
$\Box$

\

{\bf Proof of Theorem \ref{B5Apr17}.}
Since the inclusion  $\rad (R)\cap R^G \subseteq \rad (R^G)$ is obvious, it remains to show that the reverse inclusion holds, i.e., $\rad (R)\cap R^G \supseteq \rad (R^G)$. Notice that 
$\rad (\bR)=\rad (R/\rad(R))=0$. Applying Theorem \ref{A5Apr17} for the pair $\bR$, $\overline G$, gives $\rad (\bR^{\overline G})=0$. By the assumptions, $\bR^{\overline G}=\overline{R^G}$. Therefore, $$0=\rad (\overline{R^G})= \rad(R^G/\rad(R)\cap R^G)= \rad(R^G)/\rad(R)\cap R^G$$
(since $\rad(R)\cap R^G\subseteq \rad (R^G)$), i.e., $\rad (R^G)= \rad(R)\cap R^G$, as required.  $\Box$\\

{\bf Proof of Theorem \ref{8Apr17}.} Suppose that assumption 1 (resp., 2) of the theorem holds. Then, by Theorem \ref{C1.5Mont} (resp., Theorem \ref{N-2}), the ring $R^G$ is a semiprime ring.
 By Lemma \ref{c6Apr17}.(4), the ring $R^G$ is left or right  Artinian. Therefore, the ring $R^G$ is a semisimple Artinian ring. $\Box$

\

\noindent{\bf Acknowledgements.}  VB is partly supported  by Fapesp grant (2017/02946-0). VF is partly 
supported  by  CNPq grant (301320/2013-6) and by 
Fapesp grant (2014/09310-5).  This work 
 was done during the visit of the first author to the University of S\~ao Paulo whose hospitality and support are greatly acknowledged.

\small{

Department of Pure Mathematics   \;\;\;\;\;\;\;\;\;\;\;\;\;\;\;\;\;\;\;\;\;\;\;\;\;\;\;\; \;\;\;\;\;\;\;\;\;\;\;\;\;\;\;Instituto de Matem\'atica e Estat\'{i}stica

University of Sheffield  \;\;\;\;\;\;\;\;\;\;\;\;\;\;\;\;\;\;\;\;\;\;\;\;\;\;\;\;\;\;\;\;\;\;\;\;\;\;\;\;\;\;\;\;\;\;\;\; \;\;\;\;\;\;\;\;\;\;\;\;  Universidade de S\~ao Paulo

Hicks Building \;\;\;\;\;\;\;\;\;\;\;\;\;\;\;\;\;\;\;\;\;\;\;\; \;\;\;\;\;\;\;\;\;\;\;\;\;\;\;\;\;\;\;\;\;\;\;\;\;\;\;\;\;\;\;\;\;\;\;\;\;\;\;\;\; \;\;\;\;\;\,   Caixa Postal 66281

Sheffield S3 7RH \;\;\;\;\;\;\;\;\;\;\;\;\;\;\;\;\;\;\;\;\;\;\;\;\;\;\;\;\;\;\;\;\;\;\;\; \;\;\;\;\;\;\;\;\;\;\;\;\;\;\;\;\;\;\;\;\;\;\;\;\;\;\;\;\;\;\;\; S\~ao Paulo, CEP

UK    \;\;\;\;\;\;\;\;\;\;\;\;\;\;\;\;\;\;\;\;\;\;\;\;\;\;\;\;\;\;\;\;\;\;\;\; \;\;\;\;\;\;\;\;\;\;\;\;\;\;\;\;\;\;\;\;\;\;\;\;\;\;\;\;\;\;\;\;\;\;\;\;\;\;\;\;\;\;\;\;\;\;\;\;\;\;\;\;\; 05315-970, Brasil

email: v.bavula@sheffield.ac.uk \;\;\;\;\;\;\;\;\;\;\;\;\;\;\;\;\;\;\;\;\;\;\;\;\;\;\;\;\;\;\;\;\;\;\;\;\;\;\;\;\;\;\;\;\;\;\;   email: futorny@ime.usp.br    }

\end{document}